\input amsppt.sty
\magnification=\magstep1
\vsize=22.2truecm
\baselineskip=16truept
\NoBlackBoxes
\TagsOnRight
\nologo
\pageno=1
\def\N{\Bbb N}
\def\Z{\Bbb Z}

\def\l{\left}
\def\r{\right}
\def\bg{\bigg}
\def\colon{{:}\;}
\def\({\bg(}
\def\[{\bg[}
\def\){\bg)}
\def\]{\bg]}
\def\t{\text}
\def\f{\frac}
\def\mo{\roman{mod}}
\def\ord{\roman{ord}}
\def\em{\emptyset}
\def\se {\subseteq}
\def\sp {\supseteq}
\def\sm{\setminus}

\def\eq{\equiv}
\def\cs{\cdots}
\def\ls{\leqslant}
\def\gs{\geqslant}
\def\al{\alpha}

\def\da{\delta}
\def\Da{\Delta}

\def\Proof{\noindent{\it Proof}}
\def\Remark{\medskip\noindent{\it  Remark}}

\topmatter
\hbox{J. Algebra 273(2004), no. 1, 153--175.}
\bigskip
\title On the Herzog-Sch\"onheim Conjecture
for Uniform Covers of Groups\endtitle
\rightheadtext{The Herzog-Sch\"onheim Conjecture for Covers of Groups}
\author Zhi-Wei Sun\endauthor
\affil Department of Mathematics, Nanjing University
     \\Nanjing 210093, The People's Republic of China
    \\{\it E-mail:} {\tt zwsun\@nju.edu.cn}
    \\Homepage: {\tt http://pweb.nju.edu.cn/zwsun}
\endaffil
\date Communicated by Michel Brou\'e\enddate
\dedicatory Received 8 April 2002\enddedicatory
\abstract Let $G$ be any group and $a_1G_1,\ldots,a_kG_k\ (k>1)$ be left cosets in $G$.
In 1974 Herzog and Sch\"onheim conjectured that if $\Cal A=\{a_iG_i\}_{i=1}^k$
is a partition of $G$ then the (finite) indices $n_1=[G:G_1],\ldots,n_k=[G:G_k]$ cannot
be pairwise distinct. In this paper we show that if $\Cal A$ covers all the elements
of $G$ the same number of times and $G_1,\ldots,G_k$ are subnormal subgroups of $G$
not all equal to $G$,
then $M=\max_{1\ls j\ls k}|\{1\ls i\ls k\colon n_i=n_j\}|$ is not less than
the smallest prime divisor of $n_1\cs n_k$, moreover $\min_{1\ls i\ls k}
\log n_i=O(M\log^2M)$
where the $O$-constant is absolute.
\endabstract
\thanks 2000 {\it Mathematics Subject Classification.} Primary 20D60;
Secondary 05A18, 11B25, 11N45, 20D20, 20D35, 20E15, 20F16.
\newline
\indent The research was supported by
the Teaching and Research Award Program for Outstanding Young Teachers
in Higher Education Institutions of MOE, and
the National Natural Science Foundation of P. R. China.
\endthanks
\endtopmatter

\document
\heading {1. Introduction}\endheading
Let $G$ be a (multiplicative) group. As usual we use $e$ to denote
the identity element of $G$.
A left coset of a subgroup $H$ in $G$ is in the form $aH=\{ah\colon h\in H\}$
where $a\in G$.
For a finite system
$$\Cal A=\{a_iG_i\}^k_{i=1}\tag1.1$$
of left cosets in $G$, if
$$w_{\Cal A}(x)=|\{1\ls i\ls k\colon x\in a_iG_i\}|\tag1.2$$
does not depend on $x\in G$ then we call (1.1)
a {\it uniform cover} of $G$.
Only in the case $G_1=\cs=G_k=G$, (1.1) is regarded as
a {\it trivial} uniform cover of $G$.
If $ w_{\Cal A}(x)=1$
for all $x\in G$
then we call (1.1) a {\it disjoint cover} (or {\it partition}) of $G$.
A uniform cover may have no disjoint subcover (cf. [Gu]).

Any infinite cyclic group is isomorphic to the additive group $\Z$ of the integers.
The subgroups of $\Z$ different from $\{0\}$ are in the form
$n\Z=\{nx\colon x\in\Z\}$ where $n\in\Z^+=\{1,2,3,\cs\}$.
For any positive integer $n$,
the index of $n\Z$ in $\Z$ is $n$ and a coset of $n\Z$ in $\Z$ is just a residue
class
$$a+n\Z=\{x\in\Z\colon x\eq a\ (\mo\ n)\}\quad\ \ \t{where}\ a\in\Z.$$
A finite system
$$A=\{a_i+n_i\Z\}^k_{i=1}\qquad(n_1\ls\cs\ls n_k)\tag1.3$$
of residue classes is called a cover of $\Z$
if $\bigcup_{i=1}^k a_i+n_i\Z=\Z$.
Such covers were introduced by
P. Erd\H os ([E1]) in the early 1930's,
they have many surprising applications
(see, e.g. [Cr], [Gr], [Sc], [Su7], [Su9] and [Su10]).
Soon after his invention of the concept of cover of $\Z$, Erd\H os
made a conjecture that (1.3) cannot be a partition of $\Z$ if $1<n_1<\cs<n_k$. This was
confirmed by H. Davenport, L. Mirsky, D. Newman and R. Rado
(see [E2] and [NZ])
who used analysis to show that if (1.3) forms a partition of $\Z$ with $k>1$
then $n_{k-1}=n_k$.
The reader may consult [Su4], [Su5] and [Su6]
for progress on uniform covers of $\Z$.

 In the 1950's B. H. Neumann ([N1], [N2])
studied groups as unions
of cosets of subgroups while he didn't know number-theoretic research
on covers of
$\Z$. A basic result of Neumann [N1] is as follows:
If (1.1) forms a cover of a group $G$ by left cosets
but none of its proper subsystems does, then
$[G:\bigcap_{i=1}^kG_i]\ls c_k$ where $c_k$ is a constant depending on $k$.
In 1987 M. J. Tomkinson [To]
strengthened the Neumann result by showing that we can take $c_k=k!$.
By Corollary 1 of the author [Su1], for any uniform cover (1.1) of a group $G$
we also have $[G:\bigcap_{i=1}^kG_i]\ls k!$.

In 1958 S. K. Stein [St] suggested that investigations on covers of $\Z$
should be carried out on covers of abstract groups.
In 1974 M. Herzog and J. Sch\"onheim [HS]
proposed the following generalization of Erd\H os' conjecture.

\proclaim{Herzog--Sch\"onheim Conjecture}
Let $(1.1)$ be a partition of a group $G$ into $k>1$ left cosets.
Then at least two of the finite indices
$[G:G_1],\ldots,[G:G_k]$ are equal.
\endproclaim

M. M. Parmenter [Pa] and R. Brandl [Br] partially told us
when all the subgroups $G_i$ in a partition (1.1)
of group $G$ are equal or conjugate in $G$.
The Herzog-Sch\"onheim conjecture can be extended to uniform
covers of groups.

A finite group $G$ is said to be {\it pyramidal} if it contains a chain
$\{e\} =H_0\subset  H_1\subset \cs\subset  H_n=G$
of subgroups such that $[H_1:H_0]\gs\cs\gs[H_n:H_{n-1}]$
are primes in non-ascending order. In such a chain
$H_{i-1}$ is normal in $H_i$
since $[H_i:H_{i-1}]$ is the smallest prime dividing $|H_i|$
(see [Ro, 4.18]), therefore the chain
of $H$'s forms a composition series from $\{e\}$ to $G$.
Thus pyramidal groups are solvable. In 1987
M. A. Berger, A. Felzenbaum and A. S. Fraenkel [BFF4]
verified the Herzog-Sch\"onheim conjecture for pyramidal groups.

In the 1950's Erd\H os proposed the following famous unsolved problem (see [Gu]):
 Whether for any arbitrarily large $c>0$
there exists a cover (1.3) of $\Z$
 satisfying $c<n_1<\cs<n_k$?
A more general question is as follows:

\proclaim {Open Question} Let $G$ be a group and let $M$ be
 a given positive integer.
 Whether for any $N>0$ there is a finite cover $(1.1)$ of $G$ with
  each of the indices $n_i=[G:G_i]$ greater than $N$ and occurring
  at most $M$ times?
\endproclaim

For uniform covers of groups
by cosets of subnormal subgroups, we are going to confirm the generalized
Herzog-Sch\"onheim
conjecture and answer the above open question negatively! Actually we will
make further progress.

Let's introduce our basic notations.

 For $n\in\Z^+$ we let $P(n)$
be the set of prime divisors of $n$.  For a prime $p$
and a positive integer $n$, by $\ord_p n$ we mean the largest
integer $h$ such that $p^h\mid n$. For $n_1,\ldots,n_k\in\Z^+$,
$(n_1,\ldots,n_k)$ (or $(n_i)_{1\ls i\ls k}$) and $[n_1,\ldots,n_k]$ (or $[n_i]_{1\ls i\ls k}$)
stand for their greatest common divisor and least common multiple respectively.
For a real number $x$
the integral part of $x$ is denoted by $\lfloor x\rfloor$.
We also adopt conventional symbols $\sim$, $o$ and $O$ in analytic
number theory (see, e.g. [Ap]). For convenience we regard $\sum_{i\in\em}x_i$
and $\prod_{i\in\em}x_i$
as $0$ and $1$ respectively.

For a subgroup $H$ of a group $G$,
let $H_G$ denote the core (i.e. normal interior)
of $H$ in $G$, and
let $G/H$ stand for the quotient group $\{xH\colon x\in G\}$
if $H$ is normal in $G$. For a union $X$ of some left cosets of the subgroup $H$,
by $[X:H]$ we mean the number of left cosets of $H$ contained in $X$.
Sylow $p$-subgroup and Hall $\omega$-subgroup have their usual meanings
where $p$ is a prime and $\omega$ is a set of primes (cf. [Ro]). When
group $G$ and subgroups $G_1,\ldots,G_k$ are given, we let $\bigcap_{i\in I}G_i$
make sense for all $I\se\{1,\ldots,k\}$
by regarding $\bigcap_{i\in \em}G_i$ as $G$.

The main result of this paper is Theorem 4.3,
for the sake of clarity we state here a simpler version.
\proclaim{Theorem 1.1}
 Let $(1.1)$ be a nontrivial uniform cover of a group $G$
with
$$n_1=[G:G_1]\ls\cs\ls n_k=[G:G_k].\tag1.4$$
Suppose that
all the $G_i$ are subnormal in $G$, or
$G/H$ is a solvable group having a normal Sylow $p$-subgroup
where $H$ is the largest normal subgroup of $G$
contained in all the $G_i$ and  $p$ is the largest prime divisor of $|G/H|$.
Then the indices $n_1,\ldots,n_k$ cannot be pairwise
distinct. Moreover, if $|\{1\ls i\ls k\colon n_i=n\}|\ls M$ for all $n\in\Z^+$ then
we have
$$\log n_1\ls\f{e^{\gamma}}{\log 2}M\log^2M+O(M\log M\log\log M)\tag1.5$$
where the logarithm has the natural base $e=2.718...$, $\gamma=0.577...$
is the Euler constant and the $O$-constant is absolute.
\endproclaim

The next section contains some useful lemmas concerning
indices of subgroups and normal Hall subgroups.
In Sections 3 we are going to study unions of cosets.
We will investigate uniform covers
and obtain the main results in the last section.

\heading{2. Lemmas on Indices of Subgroups and Normal Hall Subgroups}\endheading

 Lemma 3.1(ii) of [Su8] can be restated as follows.

\proclaim{Lemma 2.1} Let $G$ be a group and
 $G_1,\ldots,G_k$ be subnormal subgroups of $G$ with finite index.
Then $[G:\bigcap_{i=1}^kG_i]\mid \prod_{i=1}^k [G:G_i]$
and hence
$$P\(\[G:\bigcap_{i=1}^kG_i\]\)=\bigcup_{i=1}^kP([G:G_i]).\tag2.1$$
\endproclaim

\Remark\ 2.1. If $G_1,\ldots,G_k$ are subgroups of a group $G$ with finite index,
then $[G:\bigcap_{i=1}^kG_i]\ls\prod_{i=1}^k[G:G_i]<\infty$ by Poincar\'e's theorem.
Lemma 2.1 can be viewed as an important number-theoretic property of subnormality,
it is the main reason
why covers involving subnormal subgroups are better behaved than general covers.

\proclaim{Lemma 2.2} Let $G$ be a group and
$H$ be a subnormal subgroup of $G$ with finite index. Then
$$P(|G/H_G|)=P([G:H]).\tag2.2$$
\endproclaim
\Proof. Let $\{a_iH\}^{k}_{i=1}$ be a partition of $G$ into
left cosets of $H$.
Then
$H_G=\bigcap_{g\in G}gHg^{-1}=\bigcap_{i=1}^{k}\bigcap_{h\in H}
a_ihHh^{-1}a_i^{-1}=\bigcap_{i=1}^{k}a_iHa_i^{-1}$.
Since those $a_iHa_i^{-1}$ are subnormal subgroups with index $k=[G:H]$,
(2.2) follows from Lemma 2.1. \qed

\proclaim{Corollary 2.1} Let $G$ be a finite group
and $H$ be a  Hall subgroup
of $G$. If $H$ is subnormal in $G$, then
$H$ must be normal in $G$.
\endproclaim
\Proof. By Lemma 2.2, $P(|G/H_G|)=P([G:H])$. So no prime factor of $|H|$
can divide $|G/H_G|=[G:H]|H/H_G|$. Thus $H$ coincides with $H_G$. \qed

\proclaim{Lemma 2.3} Let $G$ be a group and $H,K$ be normal
subgroups of $G$ with finite index. Let $\omega$ be a set of primes.
Then both $G/H$ and $G/K$
have normal Hall $\omega$-subgroups, if and only if $G/(H\cap K)$ has
a normal Hall $\omega$-subgroup.
\endproclaim
\Proof. Suppose that $G/(H\cap K)$ has a normal Hall $\omega$-subgroup
$F/(H\cap K)$ where $F\sp (H\cap K)$. Then $F$ is normal in $G$
and $FH/H$ is normal in $G/H$. Observe that $|FH/H|=|F/(F\cap H)|$
divides $|F/(H\cap K)|$ and hence $P(|FH/H|)\se\omega$. As $[G:FH]$ divides
$[G:F]$, $FH/H$ is a Hall $\omega$-subgroup of $G/H$. Similarly, $G/K$
has a normal Hall $\omega$-subgroup.

Now assume that $G/H$ and $G/K$ have normal Hall $\omega$-subgroups
$H^*/H$ and $K^*/K$ respectively. Then $(H^*\cap K^*)/(H\cap K)$
is normal in $G/(H\cap K)$.
Let $\bar{\omega}$ be the set of primes not in $\omega$. In light
of Lemma 2.1,
$$P([G:H^*\cap K^*])=P([G:H^*])\cup P([G:K^*])\se\bar{\omega}.$$
As $H\cap K=(H^*\cap K)\cap (H\cap K^*)$, we have
$$\align &P(|(H^*\cap K^*)/(H\cap K)|)=P([H^*\cap K^*:H^*\cap K])
\cup P([H^*\cap K^*:H\cap K^*])
\\=&P([(H^*\cap K^*)K:K])\cup P([(H^*\cap K^*)H:H])
\se P(|K^*/K|)\cup P(|H^*/H|)\se\omega.
\endalign$$
So $(H^*\cap K^*)/(H\cap K)$
is a Hall $\omega$-subgroup of $G/(H\cap K)$.
We are done. \qed

Let $G$ be a finite group and $p$ be a prime number.
Then Sylow $p$-subgroups of $G$ are just Hall $\{p\}$-subgroups of $G$.
If $G$ has a normal Sylow $p$-subgroup $S$,
then by Sylow's theorem (cf. [Ro, 5.9])
$S$ is the only Sylow $p$-subgroup of $G$.

\proclaim{Lemma 2.4} Let $G$ be a group and $H$ a subgroup of $G$ with finite
index. Then, for any $p\in P(|G/H_G|)\setminus P([G:H])$,
$G/H_G$ doesn't have a normal Sylow $p$-subgroup.
\endproclaim
\Proof. Let $K/H_G$ be a Sylow $p$-subgroup of $H/H_G$ where $K$ is a subgroup of $H$
containing $H_G$. Since $p\nmid [G:H]$, $K/H_G$ is also a Sylow $p$-subgroup
of $G/H_G$. If $G/H_G$ has a normal Sylow $p$-subgroup, then
$K/H_G$ is the unique Sylow $p$-subgroup of $G/H_G$
and therefore $K$ is normal in $G$, thus $K\se H_G$ and hence
$|K/H_G|=1$, this leads to a contradiction since $p\mid |G/H_G|$. \qed

\Remark\ 2.2. We can extend Lemma 2.4 as follows:
Let $G$ be a group and $\omega$ be a set of primes.
If $H$ is a subgroup of $G$ with finite index and
$G/H_G$ has a normal Hall $\omega$-subgroup, then
$P([G:H])\cap\omega\not=\em$ if and only if $P(|G/H_G|)\cap\omega\not=\em$.

\proclaim{Lemma 2.5} Let $G$ be a finite group and $p$ be a prime dividing $|G|$.
Then $G$ is a solvable group with a normal
Sylow $p$-subgroup, if and only if there is a composition series
$\{e\}=H_0\subset H_1\subset\cs\subset
H_n=G$ from $\{e\}$ to $G$ for which all quotients $H_1/H_0,\ldots,H_n/H_{n-1}$
have prime order, and if a
quotient is not of order $p$ then neither is the next quotient.
\endproclaim
\Proof. For the `only if' direction,
we suppose that $G$ is solvable and that $S$ is a normal
Sylow $p$-subgroup of $G$. By [Ro, 5.31] there must be a composition series
from $\{e\}$ to
$p$-group $S$ whose quotients are of order $p$. As $G/S$ is solvable
 and each prime divisor of $|G/S|$ is different from $p$,
there exists a composition series from $S$ to $G$ such that the order of
any quotient is a prime other than $p$.
Combining these we obtain a desired composition series from $\{e\}$ to $G$.

Now we consider the `if' direction.
Let $\{e\}=H_0\subset H_1\subset\cs\subset H_i\subset H_{i+1}\subset\cs\subset H_n=G$
be a composition series from $\{e\}$ to $G$ for which
$|H_1/H_0|=\cs=|H_i/H_{i-1}|=p$ and $|H_{i+1}/H_i|,
\cs,|H_n/H_{n-1}|$ are primes different from $p$. Observe that
$H_i$ is a Sylow $p$-subgroup of $G$.
By Corollary 2.1 subnormal subgroup $H_i$ is normal in $G$. So
$G$ is a solvable group with normal Sylow $p$-subgroup $H_i$.

The proof of Lemma 2.5 is now complete. \qed

\proclaim{Corollary 2.2} Let $G$ be any pyramidal group. For
the largest prime factor $p$ of $|G|$,
 $G$ has a normal Sylow $p$-subgroup.
\endproclaim
\Proof. This follows immediately from Lemma 2.5. \qed

\Remark\ 2.3. We can show that a group is pyramidal if and only if it has
a Sylow tower. Also, if a group is pyramidal then
so are its subgroups and quotient groups.

\heading {3. On Unions of Cosets}\endheading
In [Su1] it was asked whether for subgroups $G_1,\ldots,G_k$
and elements $a_1,\ldots,a_k$ of a finite
group $G$ we always have
$$\bigg|\bigcup_{i=1}^k a_iG_i\bigg|
\gs \bigg| \bigcup_{i=1}^k G_i\bigg|.$$
In 1991 Tomkinson  gave a negative answer
for $G=C_2\times C_2$ where $C_2$ is the cyclic group of order $2$.
On the other hand, we have

\proclaim{Theorem 3.1} Let $G$ be a group and $H$ its subgroup with
$[G:H]<\infty$. Let $G_1,\ldots,G_k$
be subgroups of $G$ containing $H$. Assume that either $G_1,\ldots,G_k$ are
subnormal in $G$ or there is a composition series from $H$ to $G$ whose
quotients have prime order.
Then for any $a_1,\ldots,a_k\in G$ we have
$$\aligned&|\{xH\colon x\in a_iG_i\ \t{ for some}\ i=1,\ldots,k\}|
\\\gs &|\{0\ls n<[G:H]\colon [G:G_i]\mid n\ \t{ for some}\ i=1,\ldots,k\}| ,
\endaligned$$
i.e.,
$$\[\bigcup_{i=1}^ka_iG_i:H\] \gs
\[ \bigcup_{i=1}^k[G:G_i]\Z:[G:H]\Z\].\tag3.1$$
\endproclaim

To prove it we need some preparations.

For $R\se\Z^+$  we define
$$ D(R)=\{d\in \Z^+ \colon d\mid m \ \t{for some}\ m\in R\};$$
if $k\in \Z^+$ then $R\se D(R)\se D(kR)$
where $kR=\{kr\colon r\in R\}$.
Obviously $D(\em)=\em$ and
$D(R_1\cup R_2)=D(R_1)\cup D(R_2)$ for $R_1,R_2\se\Z^+$.

Following Berger et al. [BFF4], we introduce a measure $\mu$ on
finite subsets of $\Z^+$ through
$\mu(\{ m\} )=\varphi(m)$ where $\varphi$ is  Euler's totient function.
For $m=1,2,3,\cs$ Gauss' identity $\sum_{d\mid m} \varphi(d)=m$
shows that $\mu(D(\{m\}))=m$.

For $k,m\in \Z^+$ and any finite $R\se \Z^+$ we have
$$\align\mu(D(k(R\cup \{m\})))&=\mu(D(kR)\cup D(\{km\}))
\\&=\mu(D(kR))+\mu(D(\{km\}))-\mu(D(kR)\cap D(\{km\}))
\\&=\mu(D(kR))+km-\mu(D(\{(kr,km)\colon r\in R\}))
\\&=\mu(D(kR))+k\mu(D(\{m\}))-\mu(D(kR'))
\endalign$$
where $R'=\{(r,m)\colon r\in R\}$. From this by induction we can establish
\proclaim{Lemma 3.1} Let $k$ be a positive integer and $R$ be a finite subset of $\Z^+$. Then
$$\mu(D(kR))=k\mu(D(R)).\tag3.2$$
\endproclaim

\Remark\ 3.1. The lemma was first observed by Berger et al. [BFF4].

\proclaim{Lemma 3.2} Let $\Gamma$ be a family of finite sets such that
whenever $S,T\in \Gamma$ one has $S\cap T\in \Gamma$ and $|S\cap T|=(|S|,|T|)$.
For any finite subfamily
$\Da$ of $\Gamma$ we have
$$\bg|\bigcup_{S\in \Da }S\bg|=\mu(D(\{|S|\colon S\in \Da \})).\tag3.3$$
\endproclaim
\Proof. Since $\mu(D(\em ))=\mu(\em )=0$, (3.3) holds trivially if $\Da$ is empty .

Now let $\Da_0\se\Gamma$ have cardinality $n\in\Z^+$ and assume (3.3)
for any $\Da \se \Gamma $ with
smaller cardinality. Suppose $T\in \Da_0$ and let $\Da'_0 =\Da_0 \sm\{T\}$.
By the induction hypothesis,
we have
$$\align&
 \bg|\bigcup_{S\in \Da_0}S\bg|=\bg| T\cup \bigcup_{S\in \Da_0'}S\bg|
 =\bg|\bigcup_{S\in \Da '_0}S\bg|+
|T|-\bg|\bigcup_{S\in \Da'_0}S\cap T\bg|
\\=&\mu(D(\{|S|\colon S\in \Da'_0\} ))+\mu(D(\{|T|\} ))
-\mu(D(\{|S\cap T|\colon S\in \Da'_0\}))
\\=&\mu(D(\{|S|\colon S\in \Da'_0\}))+\mu(D(\{|T|\}))
-\mu(D(\{|S|\colon S\in\Da'_0\})\cap D(\{|T|\}))
\\=&\mu(D(\{|S|\colon S\in\Da'_0\})\cup D(\{|T|\} ))=\mu(D(\{|S|\colon S\in\Da_0\})).
\endalign$$
This concludes the proof by induction. \qed

\proclaim{Lemma 3.3} Let $G$ be a group and $H$ be a subgroup of
 $G$ with finite index. Suppose that $G_1,
\cs,G_k$ are subgroups of $G$ containing $H$. Then
$$\aligned&\mu(D(\{[G_i:H]\colon 1\ls i\ls k\}))
\\=&|\{0\ls  n<[G:H]\colon [G:G_i]\mid n\  \t{for some}\ 1\ls i\ls k\}|.
\endaligned\tag3.4$$
\endproclaim
\Proof. Clearly $N=[G:H]$ is a multiple of those $n_i=[G:G_i]$ with $1\ls i\ls k$.
For each divisor $d$ of $N$ we let $X_d=\{0\ls x<N\colon x\in d\Z\}$.
If $m,n\in\Z^+$ divide $N$, then
$X_m\cap X_n=X_{[m,n]}$ has cardinality $N/[m,n]=(|X_m|,|X_n|)$.
Applying Lemma 3.2 to the family
$$\Gamma = \{X_d\colon d\in\Z^+ \ \t{and}\ d\mid N\},$$
we obtain that
$$\bg|\bigcup^k_{i=1}X_{n_i}\bg|=\mu\l(D(\{|X_{n_i}|\colon 1\ls i\ls k\})\r)
=\mu\l(D\l(\l\{\f N{n_1},\ldots,\f N{n_k}\r\}\r)\r).$$
So (3.4) holds. \qed

\medskip\noindent
{\it Proof of Theorem 3.1}. We use induction on $[G:H]$. If $G_1=\cs=G_k=G$, then
$$\bigg|\bigg\{aH\colon a\in \bigcup^k_{i=1} a_iG_i\bigg\}\bigg|
=[G:H]=\bigg|\bigg\{0\ls  n<[G:H]\colon n\in
\bigcup_{i=1}^k[G:G_i]\Z\bigg\}\bigg| .$$
Thus the case $[G:H] =1$ is trivial. So we proceed to the induction step with $[G:H]>1$
and assume that $G_j\not=G$ for some $1\ls j\ls k$.

{\it Case} 1. $G_1,\ldots,G_k$ are subnormal in $G$.
As $G_j\not=G$ there exists a proper maximal
normal subgroup $H^*$ of $G$ containing $G_j$. Observe that
each $G_i\cap H^*$ is subnormal in $H^*$
since $G_i$ is subnormal in $G$.

{\it Case} 2. There exists a  composition series from $H$ to $G$
whose quotients have prime order.
Since $H\ne G$ there is a normal subgroup $H^*$ of prime
index in $G$ for which there exists a composition series from $H$ to $H^*$
whose quotients are of prime order.

In either case, $H\se  G_i\cap  H^*\se  H^*$ and $[H^*:H]< [G:H]$.
Also, $G_iH^*$
coincides with $G$ or $H^*$.

Write $G/H^*=\{g_1H^*,\cs ,g_hH^*\}$ where $h=[G:H^*]$. Set
$$ I_s=\{1\ls  i\ls  k\colon a_iG_i\cap  g_sH^*\not=\em  \}\quad
\t{for}\ s=1,\ldots,h.$$
For each $i=1,\ldots,k$ clearly
$a_iG_i\cap g_sH^*\not= \em$ for some $1\ls s\ls h$, so
$$I_1\cup\cs\cup I_h=\{1,\ldots,k\}.\tag3.5$$

For $I=I_1\cap\cs\cap I_h$ we have
 $$\align I&=\{1\ls  i\ls  k\colon
 G_i\cap  a_i^{-1}g_sH^*\not=\em\ \t{for all}\ s=1,\ldots,h\}\\
&=\{1\ls  i\ls  k\colon xG_i\cap yH^*=x(G_i\cap  x^{-1}yH^*)\not=\em
\ \t{for all}\ x,y\in G\}
\\&=\{1\ls  i\ls  k\colon G_iH^*=G\}\ \ \ \t{(by Lemma 2.1 of [Su8])}.
\endalign$$
 Let $R=\{[G_i\cap  H^*:H]\colon i\in  I\}$. Then  $hR=\{[G_i:H]\colon i\in  I\}$
 since $[G_i:G_i\cap H^*]=[G_iH^*:H^*]=h$ for all $i\in I$.

Let $s\in\{1,\ldots,h\}$ and $R_s=\{[G_i \cap H^*:H]\colon i\in I_s\sm I\}$.
For $i \in I_s\sm I$, as $G_iH^*=H^*$ we have  $G_i \se
H^*$ and $a_iG_i \se g_sH^*$. So $R_s=\{[G_i:H]\colon i\in  I_s\sm I\}$.

If $i\in I_s$, then $g_s^{-1}a_iG_i \cap H^*$ is nonempty and hence it is
a left coset of $G_i\cap  H^*$ in $H^*$. Clearly
$$\bigcup^k_{i=1}a_iG_i=\bigcup_{s=1}^h\bigcup^k_{i=1}a_iG_i\cap g_sH^*
=\bigcup_{s=1}^h\bigcup_{i\in I_s}a_iG_i\cap g_sH^*$$
and so
$$\[\bigcup^k_{i=1}a_iG_i:H\]
=\sum_{s=1}^h\[\bigcup_{i\in I_s}a_iG_i\cap g_sH^*:H\]
=\sum_{s=1}^h\[\bigcup_{i\in I_s}g_s^{-1}a_iG_i\cap H^*:H\].$$
Thus, by the induction hypothesis,
$$\align  \[\bigcup_{i=1}^k a_iG_i:H\]\gs
&\sum_{s=1}^{h}\bigg|\bg\{0\ls  n< [H^*:H]\colon
n\in \bigcup_{i\in  I_s}[H^*:G_i\cap  H^*]\Z\bg\}\bigg|
\\=&\sum_{s=1}^{h} \mu(D(\{[G_i \cap H^*:H]\colon i\in  I_s\}))\ \ \t{(by Lemma 3.3)}
\\=&\sum_{s=1}^{h} \mu(D(R\cup  R_s))=h\mu(D(R))+\sum_{s=1}^{h} \mu(D(R_s)\sm D(R))
\\\gs & \mu(D(hR))+\sum_{s=1}^{  h}\mu( D(R_s)\sm D(hR)) \ \ \ \t{(by Lemma 3.1)}.
\endalign$$
It follows that
$$\align  \[\bigcup_{i=1}^k a_iG_i:H\]
\gs&\mu\(D(hR)\cup\bigcup_{s=1}^{h}\l(D(R_s)\sm D(hR)\r)\)
=\mu\(\bigcup^h_{s=1}D(hR\cup R_s)\)
\\=&\mu\(\bigcup^h_{s=1}D(\{[G_i:H]\colon i\in I_s\})\)
=\mu(D(\{[G_i:H]\colon 1\ls i\ls k\}))
\\=&\bg|\bg\{0\ls n<[G:H]\colon n\in \bigcup_{i=1}^{k}[G:G_i]\Z\bg\}\bg|
\ \ \t{(by Lemma 3.3)}.
\endalign$$
This completes the proof.   \qed

\Remark\ 3.2. A theorem of C. A. Rogers (cf. [HR]) indicates that
if $a_i\in \Z$ and $n_i\in \Z^+$
for $i=1,\ldots,k$ then for any positive multiple $N$ of $n_1,\ldots,n_k$ we have
$$\bigg|\bigg\{0\ls  x< N\colon x\in \bigcup_{i=1}^{k}a_i+n_i\Z\bigg\}\bigg|
\gs \bigg|\bigg\{0\ls  x< N\colon x\in \bigcup^k_{i=1} n_i\Z\bigg\}\bigg|,\tag3.6$$
this is just our Theorem 3.1 in the case where $G$ is the infinite cyclic group $\Z$.
(It should be mentioned that Simpson [Si] presented this as his Lemma 2.3
but gave a wrong proof.)
In view of Lemma 3.3, Lemma IV of [BFF4]
is equivalent to our Theorem 3.1 in
the case where $G$ is a pyramidal group and $H$ is the smallest subgroup $\{e\}$.

As in additive number theory, for any $S\se \Z$ we let $d(S)$ denote the asymptotic
density
$$\lim_{N\to +\infty}\f{|\{0\ls n<N\colon n\in S\}|}N$$
if the limit exists. It is easy to see that for system (1.3) we have
$$d\(\bigcup^k_{i=1}a_i+n_i\Z\)
=\f1N \bigg|\bigg\{0\ls  x< N\colon x\in \bigcup_{i=1}^{k}a_i+n_i\Z\bigg\}\bigg|$$
where $N$ is any positive multiple of $[n_1,\ldots,n_k]$.

Here we restate Lemma 2 of [Su2] (proved by
the inclusion-exclusion principle).

\proclaim{Lemma 3.4} Let $n_1,\ldots,n_k$ be positive integers
and let $P$ be a finite set of primes such that
$P(n_i)\se P$ for all $i=1,\ldots,k$. Then
$$d\(\bigcup^k_{i=1}n_i\Z\)=\(\prod_{p\in P}\f{p-1}p\)
\sum\Sb n\in\bigcup^k_{i=1}n_i\Z^+\\P(n)\se P\endSb\f1n.\tag3.7$$
\endproclaim

Now we are able to give
\proclaim{Theorem 3.2} Let $G$ be a group and $G_1,\ldots,G_k, H$ be
subgroups of $G$ with finite index. Let
$a_1,\ldots,a_k\in G$ and assume that the union of $a_1G_1,\ldots,a_kG_k$
coincides with a union of some
left cosets of $H$. Let $h=[G:H]$ and $n_i=[G:G_i]$ for $i=1,\ldots,k$. Then we have
$$\f{(n_1,\ldots,n_k)}{(h,n_1,\ldots,n_k)}
\ls\sup_{n\in  \Z^+}|\{1\ls  i\ls  k\colon n_i=n\}|
\sum_{d\mid \f{ [ n_1,\ldots,n_k]}{(n_1,\ldots,n_k)}}\f1d\tag3.8$$
in the following four cases.

{\rm (a)} All the $G_i$ are subnormal and $H$ is normal in $G$.

{\rm (b)} All the $G_i$ are normal and $H$ is subnormal in $G$.

{\rm (c)} All the $G_i$ are normal in $G$ and $G/ \bigcap^k_{i=1}G_i$ is solvable.

{\rm (d)} $H$ is normal in $G$, and $G/H$ or each $G/(G_i)_G$ is solvable.
\endproclaim
\Proof. In either case $G_iH=HG_i$ for $i=1,\ldots,k$.
Clearly $G_iH$ is subnormal in $G$ for every
$i=1,\ldots,k$ in case $(a)$ or case $(b)$, and there is a composition series
from $H$ or $\bigcap^k_{i=1}(G_i)_G$ to $G$ whose quotients have prime order
in case (c) or case (d).
Note that
$$\bigcup^k_{i=1}a_iG_i=\bigg(\bigcup^k_{i=1}a_iG_i\bigg)H =\bigcup^k_{i=1}a_iG_iH.$$

With the help of Theorem 3.1, for a suitable $F\in\{\bigcap^k_{i=1}(G_i)_G,H\}$ we have
$$\align&\sum^k_{i=1}\f1{n_i}=\sum^k_{i=1}\f{[a_iG_i:F]}{[G:F]}
\\\gs&\f1{[G:F]}\[\bigcup_{i=1}^ka_iG_i:F\]=
 \f1{[G:F]}\[\bigcup_{i=1}^ka_iG_iH:F\]
\\\gs&\f1{[G:F]}\bigg |\bigg\{0\ls  n<[G:F]
\colon [G:G_iH]\mid n\ \t{for some}\ i=1,\ldots,k\bigg\}\bigg|
\\=&d\bigg(\bigcup_{i=1}^k[G:G_iH]\Z\bigg).\endalign$$

Let $S=\{n_1,\ldots,n_k\}$, $P=\bigcup_{ n\in  S}P(n)$ and
$\bar P=\{n\in\Z^+\colon P(n)\se P\}$.
Then
$$ \bigcup_{i=1}^kP([G:G_iH])\se \bigcup_{i=1}^kP(n_i)=P,\ \
\t{i.e.}\ \ \{[G:G_iH]\colon 1\ls i\ls k\}\se\bar P.$$
Obviously $(n_1,\ldots,n_k)/(h,n_1,\ldots,n_k)=[h,(n_1,\ldots,n_k)]/h$
divides $n_i/(h,n_i)=[h,n_i]/h$ and
$[G:G_iH]$ divides $(h,n_i)$, therefore
$$\f{ n_i}{(n_1,\ldots,n_k)/(h,n_1,\ldots,n_k)}\in\bar P\cap \bigcup_{j=1}^k[G:G_jH]\Z.$$

Clearly $[ n_1,\ldots,n_k]/(n_1,\ldots,n_k)$ can be written
in the form $\prod_{p\in P}p^{\da_p}$ where
$\da_p\in\N=\{0,1,2,\cs\}$. For any $p\in P$ and $1\ls i,j\ls k$ we have
$$\ord_p n_i-\ord_p n_j
\ls \ord_p[n_1,\ldots,n_k]-\ord_p(n_1,\ldots,n_k)=\da_p.$$
So, if $n,n'\in  S$, $k_p,l_p\in\N$ and
$$n\prod_{p\in  P}p^{k_ p (1+\da_p)}=n'\prod_{p\in P} p^{ l_p (1+\da_ p)},$$
then $k_ p=l_ p$ for all $p\in  P$ and hence $n=n'.$

Let $$M=\sup_{n\in\Z^+}|\{1\ls  i\ls  k \colon n_i=n\}|
=\max_{ n\in  S}\sum^k\Sb i=1\\n_i=n\endSb 1.$$
In view of Lemma 3.4 and the above,
$$\aligned\sum^k_{i=1}\f1{ n_i}\gs & d\bigg(\bigcup_{i=1}^k [G:G_iH]\Z\bigg)
= \bigg(\prod_{ p\in  P}\f{p-1}p\bigg)\sum_{m\in\bar P\cap \bigcup_{i=1}^k[G:G_iH]\Z}\f1m
\\\gs& \bigg (\prod_{ p\in  P} \f{p-1}p\bigg)
\sum_{n\in  S}\l(\f n{(n_1,\ldots,n_k)/(h,n_1,\ldots,n_k)}\r)^{-1}
\prod_{p\in  P}\sum^{\infty}_{i=0}\f1{p^{i(1+\da_p)}}
\\=&\sum_{n\in S}\f1n\cdot\f{(n_1,\ldots,n_k)}{(h,n_1,\ldots,n_k)}
\prod_{p\in  P}\(\f{p-1}p\bigg/\bigg(1-\f1{p^{1+\da_p}}\bigg)\bigg)
\\=&\f1M\sum_{n\in S}\f Mn\cdot\f{(n_1,\ldots,n_k)}{(h,n_1,\ldots,n_k)}
\prod_{p\in P}\f{ p^{\da_p}}{1+p+\cs+p^{\da_p}}
\\\gs& \f1M\sum^k_{i=1}\f1{n_i}\cdot\f{(n_1,\ldots,n_k)}{(h,n_1,\ldots,n_k)}
\prod_{p\in P}\l(1+\f1p+\cs+\f1{p^{\da_p}}\r)^{-1}.
\endaligned$$
Therefore
$$\f{(n_1,\ldots,n_k)}{(h,n_1,\ldots,n_k)}\ls M\prod_{p\in P}\l (1+\f1p+\cs+\f1{p^{\da_p}}\r)
=M\sum_{d\mid \f{ [ n_1,\ldots,n_k]}{(n_1,\ldots,n_k)}}\f1d.$$
We are done. \qed

Our Theorem 3.2 is powerful, it will be applied in Section 4.

\heading {4. On Uniform Covers}\endheading

\proclaim{Lemma 4.1} Let $(1.1)$ be a finite system of left cosets
in a group $G$. Then
$$K_{\Cal A}=\{x\in G\colon w_{\Cal A}(gx)=w_{\Cal A}(g)\ \t{for all}\ g\in G\}\tag4.1$$
is a subgroup of $G$ containing $\bigcap_{i=1}^k G_i$.
For any nonempty subset $I$ of $\{1,\ldots,k\}$,
the union $\bigcup_{i\in I}a_iG_i$ coincides with a union of
some left cosets of $K_{\Cal A}\cap\bigcap_{j\in \bar I}G_j$
where $\bar I=\{1,\ldots,k\}\sm I$.
\endproclaim
\Proof. If $x,y\in K_{\Cal A}$ then
$w_{\Cal A}(gxy^{-1})=w_{\Cal A}(gxy^{-1}y)=w_{\Cal A}(gx)=w_{\Cal A}(g)$
for all $g\in G$. So $K_{\Cal A}$ is a subgroup of $G$.
For $g\in G$ and $x\in G_i$, clearly
$gx\in a_iG_i$ if and only if $g\in a_iG_ix^{-1}=a_iG_i$.
Thus $K_{\Cal A}\sp\bigcap_{i=1}^k G_i$.

Let $g\in G$ and $x\in K_{\Cal A}\cap\bigcap_{j\in\bar I}G_j$.
For $j\in \bar I$ we have
$gx\in a_jG_j\Leftrightarrow g\in a_jG_j$. Therefore
$$\align&|\{i\in I\colon gx\in a_iG_i\}|=
w_{\Cal A}(gx)-|\{j\in\bar I\colon gx\in a_jG_j\}|
\\=& w_{\Cal A}(g)-|\{j\in\bar I\colon g\in a_jG_j\}|
=|\{i\in I\colon g\in a_iG_i\}|.\endalign$$
It follows that
$$g\in\bigcup_{i\in I}a_iG_i\Longrightarrow
g\(K_{\Cal A}\cap\bigcap_{j\in\bar I}G_j\)\se\bigcup_{i\in I}a_iG_i.$$
So $X=\bigcup_{i\in I}a_iG_i$ is identical with
$\bigcup_{g\in X}g(K_{\Cal A}\cap\bigcap_{j\in \bar I}G_j)$. We are done. \qed

\proclaim{Theorem 4.1} Let $(1.1)$ be a nontrivial uniform cover of
a group $G$ by left cosets.
Let $n_i=[G:G_i]$ for $i=1,\ldots,k$ and $[n_1,\ldots,n_k]=\prod^r_{t=1}p_t^{\al_t}$
where $p_1,\ldots,p_r$
are distinct primes and $\al_1,\ldots,\al_r$ are positive integers. Let
$$\beta_r=\min\{1\ls \beta\ls\al_r\colon \beta=\ord_{p_r}n_i
\ \t{for some}\ i=1,\ldots,k\},\tag4.2$$
$$\varepsilon_r=\(1-\f1{p_r^{\al_r-\beta_r+1}}\)
\prod_{0<t<r}\(1-\f1{p_t^{\al_t+1}}\)\tag4.3$$
and $$M_r=\max\l\{|\{1\ls i\ls k\colon n_i=n_j\}|
\colon 1\ls j\ls k\ \&\ p_r\mid n_j\r\}.\tag4.4$$
Then we have
$$p^{\beta_r}_r\ls \varepsilon_r M_r\prod^{r}_{t=1}\f{p_t}{p_t-1}\tag4.5$$
providing the following $(a)$ and $(b)$, or $(c)$ in the case $p_1<\cs<p_r$.

{\rm (a)} If not all the $G_i$ with $p_r\mid n_i$ are subnormal in $G$,
then all the $G/(G_i)_G$ with $p_r\mid n_i$,
or those with $p_r\nmid n_i$, are solvable.

{\rm (b)} For each $i$ with $n_i>p_r$ and $p_r\nmid n_i$, if $G_i$ is not subnormal in $G$
then $G/(G_i)_G$
has a normal Sylow $p_r$-subgroup.

{\rm (c)} {\it $\bar G=G/(\bigcap^k_{i=1}G_i)_G$ is a solvable group
having a normal Sylow $p$-subgroup where
$p$ is the largest prime divisor of $|\bar G|$.}
\endproclaim
\Proof. Suppose that $p_1<\cs<p_r$ and (c) holds. Since
$(\bigcap_{i=1}^kG_i)_G=\bigcap_{i=1}^k(G_i)_G$, by Lemma 2.3
each $G/(G_i)_G$ is a
solvable group having a normal Sylow $p$-subgroup.
In view of Lemma 2.4, if  $p\mid|G/(G_i)_G|$ then $p\mid n_i$. On the other hand
$$ P(|\bar G|)=\bigcup^k_{i=1} P(|G/(G_i)_G|)
\supseteq \bigcup_{i=1}^k P(n_i)=\{p_1,\ldots,p_r\},$$
so we have $p=p_r$. Therefore both (a) and (b) hold.

Below we prove (4.5) under the conditions (a) and (b).

Let $I=\{1\ls i\ls k\colon p_r\mid n_i\}$ and $\bar I=\{1,\ldots,k\}\sm I$. Since
$ w_{\Cal A}$ is constant, $K_{\Cal A}=G$. By Lemma 4.1, $\bigcup_{i\in I}a_iG_i$ coincides
with a union of
some left cosets of $\bigcap_{j\in \bar I}G_j$.
Let $H=(\bigcap_{j\in\bar I}G_j)_G=\bigcap_{j\in\bar I}(G_j)_G$.
Then $G/H$ is finite and $\bigcup_{i\in I}a_iG_i$ is a union of finitely many cosets of $H$.
Note that $G/H$ is solvable if and only if $G/(G_j)_G$ is solvable for all $j\in\bar I$
(cf. [Ro, 7.46 and 7.50]). By condition $(a)$ and Theorem 3.2, we have
$$\f{(n_i)_{i\in I}}{(|G/H|,(n_i)_{i\in I})}\ls
\sup_{n\in\Z^+}|\{i\in I\colon n_i=n\}|\sum_{d\mid\f{[n_i]_{i\in I}}
{(n_i)_{i\in I}}}\f1d.$$
Therefore
$$\aligned\f{[|G/H|,p_r^{\beta_r}]}{|G/H|}\ls&\f{[|G/H|,(n_i)_{i\in I}]}{|G/H|}
\ls M_r\sum_{d\mid
p_r^{\al_r-\beta_r}\prod_{0<t<r}p_t^{\al_t}}\f1d
\\=&M_r\(1+\f1{p_r}+\cs+\f1{p_r^{\al_r-\beta_r}}\)
\prod_{0<t<r}\(1+\f1{p_t}+\cs+\f1{p_t^{\al_t}}\)
\\=&M_r\f{p_r^{\al_r-\beta_r+1}-1}{p_r^{\al_r-\beta_r}(p_r-1)}
\prod_{0<t<r}\f{p_t^{\al_t+1}-1}{p_t^{\al_t}(p_t-1)}
=\varepsilon_r M_r\prod^r_{t=1}\f{p_t}{p_t-1}.\endaligned$$

Now it suffices to show that $p_r\nmid|G/H|$ under condition $(b)$.
In view of Lemma 2.1, $P(|G/H|)=\bigcup_{j\in\bar I}P(|G/(G_j)_G|)$. Let $j\in \bar I$.
By [Ro, 4.14], $G/(G_j)_G$ can be embedded in the symmetric group
$S_{n_j}$. If $n_j<p_r$, then $|S_{n_j}|=n_j!\not\equiv0\ (\mo \ p_r)$
and hence $p_r\nmid|G/(G_j)_G|$. If $G_j$ is subnormal in $G$,
 then by Lemma 2.2, $P(|G/(G_j)_G|)=P(n_j)$
doesn't contain $p_r$.
When $G/(G_j)_G$ has a normal Sylow $p_r$-subgroup,
$p_r\nmid|G/(G_j)_G|$ by Lemma 2.4.
So $p_r\not\in P([G:H])$ and we are done. \qed

\Remark\ 4.1.
In [BS] N. Burshtein and Sch\"onheim investigated
disjoint covers of $\Z$ having moduli
occurring at most twice. In 1976 Burshtein [Bu]
conjectured that
for any disjoint cover (1.3) with
each modulus occurring at most $M\in \Z^+$ times,
if $p_1<\cs<p_r$ are the distinct prime divisors
of $[n_1,\ldots,n_k]$ then
$$p_r\ls M\prod^r_{t=1}\f{p_t}{p_t-1};\tag4.6$$
he also realized that
the smallest modulus in such a disjoint cover
cannot be arbitrarily large by his conjecture.
The conjecture was later proved by
Simpson [Si], and by Berger et al. [BFF3] independently.
In [BFF2] and [BFF4] Berger et al. obtained
the analogy of the Burshtein conjecture for partitions of
finite nilpotent groups and pyramidal groups, their results follow from
our Theorem 4.1 in view of Corollary 2.2.

\proclaim{Corollary 4.1} Let $G$ be a group of squarefree order, and
$(1.1)$ be a nontrivial uniform cover
of $G$ with $p_1<\cs<p_r$ being the prime divisors of the indices $[G:G_1],\ldots,[G:G_k]$.
Then for some
$n\equiv0\ (\mo \ p_r)$ we have
$$|\{1\ls i\ls k\colon [G:G_i]=n\}|\gs \f{p_1\cdots p_r}{\prod_{0<t<r}(p_t+1)}
\gs\max\l\{p_1,\f{2p_r}{r+1}\r\}.\tag4.7$$
\endproclaim
\Proof. As $\bar G=G/(\bigcap^k_{i=1}G_i)_G$ has squarefree
order, by [Ro, Exercise 609] $\bar G$ is a solvable group having
a normal Sylow $p$-subgroup where
$p$ is the largest prime divisor of $|\bar G|$.
By Theorem 4.1, for some $j=1,\ldots,k$ with $[G:G_j]\in p_r\Z$,
we have
$$p_r\ls M\(1-\f1{p_1^2}\)\cs\(1-\f1{p^2_{r-1}}\)\(1-\f1{p_r}\)\prod^r_{t=1}\f{p_t}{p_t-1}$$
where  $M=|\{1\ls i\ls k\colon [G:G_i]=[G:G_j]\}|$. Thus
$$\align\f M{p_r}\gs&\(1+\f1{p_1}\)^{-1}\cs\(1+\f1{p_{r-1}}\)^{-1}=
\prod_{0<t<r}\f{p_t}{p_t+1}=\prod_{0<t<r}\(1-\f1{p_t+1}\)
\\\gs&\prod_{0<t<r}\l(1-\f1{p_1+t}\r)=\f{p_1}{p_1+r-1}\gs\max\l\{\f{p_1}{p_r},\f{2}{r+1}\r\}
\endalign$$
and the desired result follows. \qed

Our progress on the Herzog-Sch\"onheim conjecture is as follows.

\proclaim{Corollary 4.2} Let $G$ be a group and $(1.1)$ be
a nontrivial uniform cover of $G$ by left cosets.
Let $r$ be the number of distinct prime divisors of
$N=[[G:G_1],\ldots,[G:G_k]]$,
and let $p$ be any prime divisor of $|G/(\bigcap^k_{i=1}G_i)_G|$ greater than $r$
(e.g. the largest prime divisor of $N$).
Suppose that all those $G_i$ with $[G:G_i]\gs p$
are subnormal in $G$ and $p$ divides $N$,
or $G/(\bigcap^k_{i=1}G_i)_G$ is a solvable group having a normal Sylow $p$-subgroup.
Then there is a pair $\{i,j\}$ with $1\ls i<j\ls k$ such that
$[G:G_i]=[G:G_j]\equiv 0\ (\mo \ p)$.
\endproclaim
\Proof.
If $\bar G=G/\bigcap^k_{i=1}(G_i)_G$  is a solvable group having
 a normal Sylow $p$-subgroup,
then so is each $G/(G_i)_G$ by Lemma 2.3, also $p$ divides $N$
by Lemmas 2.1 and 2.4.

Set $p_r=p$ and let $p_1,\ldots,p_{r-1}$ be the other $r-1$ distinct
prime divisors of $N$. By Theorem 4.1 we have
$$p_r<|\{1\ls i\ls k\colon [G:G_i]=[G:G_j]\}|\prod^r_{t=1}\f{p_t}{p_t-1}$$
for some $j\in\{1,\ldots,k\}$ with $p\mid[G:G_j]$. Thus
$$|\{1\ls i\ls k\colon[G:G_i]=[G:G_j]\}|
>p\prod^r_{t=1}\(1-\f1{p_t}\)\gs(r+1)\prod^r_{t=1}\(1-\f1{t+1}\)=1.$$
Therefore $[G:G_i]=[G:G_j]\eq0\ (\mo\ p)$ for some $i=1,\ldots,k$
with $i\ne j$. \qed

For cyclic groups we can say something more general than Theorem 4.1.
\proclaim{Theorem 4.2} Let $(1.1)$ be a nontrivial uniform cover of a cyclic
group $G$ by cosets of subgroups $G_i$ of indices $n_i$. Assume that
$[n_1,\ldots,n_k]=p^{\al_1}_1\cs p^{\al_r}_r$ where $p_1,\ldots,p_r$
are distinct primes and $\al_1,\ldots,\al_r$ are positive integers.
Let $\al$ be a positive integer in $\Lambda=\{\ord_{p_r}n_i\colon 1\ls i\ls k\}$
and $\beta$ be the largest integer in $\Lambda\cup\{0\}$
less than $\al$. Then
$$p_r^{\al-\beta}\ls\varepsilon\max\Sb 1\ls j\ls k\\p_r^{\al}\mid n_j\endSb
|\{1\ls i\ls k\colon n_i=n_j\}|\prod^r_{t=1}\f{p_t}{p_t-1}
\tag4.8$$
where
$$\varepsilon=\(1-\f1{p_1^{\al_1+1}}\)\cs\(1-\f1{p_{r-1}^{\al_{r-1}+1}}\)
\(1-\f1{p_r^{\al_r-\al+1}}\).\tag4.9$$
Consequently,
$$\max\Sb 1\ls j\ls k\\\ord_{p_r}n_j=\al_r\endSb|\{1\ls i\ls k\colon n_i=n_j\}|
\gs p_r\prod_{0<t<r}\f{p_t-1}{p_t}\gs\f{p_r}r.\tag4.10$$
\endproclaim
\Proof. Let $I=\{1\ls i\ls k\colon p_r^{\al}\mid n_i\}$ and
$\bar I=\{1,\ldots,k\}\sm I$. Set $H=\bigcap_{j\in\bar I}G_j$
and
$$M=\sup_{n\in\Z^+}|\{i\in I\colon n_i=n\}|=\max\Sb 1\ls j\ls k\\p_r^{\al}\mid n_j\endSb
|\{1\ls i\ls k\colon n_i=n_j\}|.$$
As in the proof of Theorem 4.1 we have
$$\f{[|G/H|,p^{\al}_r]}{|G/H|}\ls \varepsilon M\prod^{r}_{t=1}\f{p_t}{p_t-1}.$$
If $\bar I=\em$ then $\ord_{p_r}|G/H|=0=\beta$.
When $\bar I\ne\em$ and $G=\langle a\rangle$,
$H=\bigcap_{j\in\bar I}\langle a^{n_j}\rangle=
\langle a^{[n_j]_{j\in\bar I}}\rangle$ and therefore
$\ord_{p_r}|G/H|=\ord_{p_r}[n_j]_{j\in\bar I}=\beta$.
So $[|G/H|,p_r^{\al}]/|G/H|=p_r^{\al-\beta}$ and hence (4.8) holds.
If we take $\al=\al_r$ then
$\varepsilon\ls 1-p_r^{-1}=(p_r-1)/p_r$ and the first inequality in (4.10) follows.
For the second inequality
in (4.10), we note that
$\prod_{0<t<r}(p_t-1)/p_t\gs\prod_{0<s<r}s/(s+1)=1/r$.
This ends our proof.  \qed

\Remark\ 4.2. Let (1.3) be a disjoint cover of $\Z$
with each modulus occurring at most $M$ times.
Suppose that $p_1,\ldots,p_r$ are the distinct prime divisors of
$n_1,\ldots,n_k$.
In 1986 Simpson [Si] showed the
inequality $p_r\ls M\prod_{0<t<r}p_t/(p_t-1)$.
(In the case $r\gs2$ and $p_1<\cs<p_r$, the weaker inequality
$M\gs p_2(p_1-1)/p_1$ was first noted in [BFF1].)
This improvement to the original
Burshtein conjecture was strengthened in [Su2]
where the author got Theorem 4.2 for disjoint covers of $\Z$.
For any cyclic group $G$, Theorem 4.1 corresponds to Theorem 4.2
in the case $\al=\min (\Lambda\cap\Z^+)$.

From now on variable $p$ will only take prime values as in number
theory.

\proclaim{Lemma 4.2} For $M\gs2$, if $q>1$ is an integer with
$q<M\prod_{p\ls q}p/(p-1)$ then
$$q<e^{\gamma}M\log M+O(M\log\log M)
\ \t{and}\  \pi(q)\ls e^{\gamma}M+O(M/\log M)\tag4.11$$
where $\pi(q)$ is the number of primes not exceeding $q$
and the $O$-constants are absolute.
\endproclaim
\Proof. A well-known theorem of Mertens (see Theorem 13.13 of [Ap]) asserts that
$$\prod_{p\ls x}\(1-\f1p\)
=\f{e^{-\gamma}}{\log x}+O\(\f1{\log^2x}\)\ \quad \t{for}\ x\gs2.$$
Thus for $x\in[2,+\infty)$ we have
$$\prod_{p\ls x}\f p{p-1}=\f{ e^{\gamma}\log x}{1+O(\f1{\log x})}
=(e^{\gamma}\log x)\(1+O\(\f1{\log x}\)\)=e^{\gamma}\log x+O(1).$$
(Note that $(1-z)^{-1}=1+z/(1-z)=1+O(z)$ when $|z|<1/2$.)
For $M\gs2$ we let $c(M)$ be the smallest positive integer $x$
such that $\prod_{p\ls x}p/(p-1)\ls x/M$,
obviously $c(M)>2$.

When $M'\gs M$, we have $c(M')\gs c(M)$ because
$$\f1{c(M')}\prod_{p\ls c(M')}\f p{p-1}\ls\f1{M'}\ls\f1M.$$
If $c(M')>c(M)$ for no $M'>M$, then
$$\f1{c(M)}\prod_{p\ls c(M)}\f p{p-1}
=\f1{c(M')}\prod_{p\ls c(M')}\f p{p-1}\ls\f1{M'}\quad\t{for all}\ M'>M,$$
and hence $c(M)^{-1}\prod_{p\ls c(M)}p/(p-1)=0$ which is impossible.
So $c(M)\to+\infty$ as $M\to+\infty$.
By the definition of $c(M)$,
$$\f1{c(M)}\prod_{p\ls c(M)}\f p{p-1}\ls\f1M<\f1{c(M)-1}\prod_{p\ls c(M)-1}\f p{p-1}.$$
Thus $c(M)$ cannot be a prime, and
$$1-\f1{c(M)}=\f{c(M)-1}{c(M)}<\f M{c(M)}\prod_{p\ls c(M)-1}\f p{p-1}
=\f M{c(M)}\prod_{p\ls c(M)}\f p{p-1}\ls 1.$$
Since $\prod_{p\ls c(M)}(1-p^{-1})\to 0$ as $M\to+\infty$, we have $M=o(c(M))$.

By the above,
$$\align&e^{\gamma}M\f{\log c(M)}{c(M)}=\f M{c(M)}\prod_{p\ls c(M)}\f p{p-1}+O\(\f M{c(M)}\)
\\=&1+O\(\f1{c(M)}\)+O\(\f M{c(M)}\)= 1+O\(\f M{c(M)}\)\endalign$$
and hence
$$\f{c(M)}{\log c(M)}
=e^{\gamma}M\(1+O\l(\f M{c(M)}\r)\)^{-1}=e^{\gamma}M+O\(\f{M^2}{c(M)}\).$$
It follows that
$$\log c(M)\sim \log\f{c(M)}{\log c(M)}=\log(e^{\gamma}M(1+o(1)))\sim\log M$$
and
$$c(M)\sim e^{\gamma}M\log c(M)\sim e^{\gamma}M\log M.$$
Thus
$$\align\log c(M)=&\log(e^{\gamma}M\log M)+\log(c(M)/(e^{\gamma}M\log M))
\\=&\gamma+\log M+\log\log M+o(1)=\log M+O(\log\log M)\endalign$$
and
$$\align c(M)=&e^{\gamma}M\log c(M)+O\(\f{M^2\log c(M)}{c(M)}\)
\\=&e^{\gamma}M(\log M+O(\log\log M))+O\(\f{M^2}{e^{\gamma}M}\)
\\=&e^{\gamma}M\log M+O(M\log\log M).\endalign$$

The famous prime number theorem (see Chapter 4 of [Ap]) states that
$$\pi(x)=\sum_{p\ls x}1\sim\f x{\log x}\quad\t{as}\ x\to+\infty,$$
moreover $\pi(x)=x/\log x+O( x/\log^2x)$ (for $x\gs2$) by [Bo] or [DV].
Hence
$$\align&\pi(c(M))=c(M)/\log c(M)+O(c(M)/\log^2c(M))
\\=&e^{\gamma}M+O\(\f{M^2}{e^{\gamma}M\log M}\)
+O\(\f{e^{\gamma}M\log M}{\log^2 M}\)=e^{\gamma}M+O\(\f M{\log M}\).
\endalign$$

It is easy to see that
$$\f1{l+1}\prod_{p\ls l+1}\f p{p-1}\ls\f1l\prod_{p\ls l}\f p{p-1}
\quad\ \t{for every}\ l=1,2,3,\cs.$$
Therefore
$$\f1n\prod_{p\ls n}\f p{p-1}\ls\f1{c(M)}\prod_{p\ls c(M)}\f p{p-1}\ls\f1M$$
for all $n=c(M),c(M)+1,\cs$. When an integer $q>1$ satisfies the inequality
 $q<M\prod_{p\ls q}p/(p-1)$ (i.e. $q^{-1}\prod_{p\ls q}p/(p-1)>M^{-1}$),
we must have $q<c(M)$ and  $\pi(q)\ls\pi(c(M))$,
so (4.11) follows. This completes the proof. \qed

\proclaim{Theorem 4.3} Let $(1.1)$ be a nontrivial uniform
cover of a group $G$ such that among
the indices $[G:G_1]\ls\cs\ls[G:G_k]$ each occurs at most $M\in\Z^+$ times.
Let $p_*$ and $p^*$ be the smallest and the largest prime divisors of
$N=[[G:G_1],\ldots,[G:G_k]]$ respectively. Suppose that all the $G_i$
with $[G:G_i]\gs p^*$ are subnormal in $G$,
or $G/H$ is a solvable group having
a normal Sylow $p'$-subgroup where
$H$ is the largest normal subgroup of $G$ contained in all the $G_i$ and
$p'$ is the greatest prime divisor of $|G/H|$ (equivalently,
 there is a composition series from $H=(\bigcap^k_{i=1}G_i)_G$ to $G$
whose quotients have prime order and if
a quotient is not of the maximal order then neither is the next quotient).
Then we have the following {\rm (i)--(iv)}
with the $O$-constants absolute.

{\rm (i)} $M\gs p_*$, moreover among the $k$ indices
$[G:G_1],\ldots,[G:G_k]$ there exists a multiple of $p^*$ occurring at least
$1+\lfloor p^*\prod_{p\mid N}(p-1)/p\rfloor\gs p_*$ times.

{\rm (ii)} All prime divisors of $[G:G_1],\ldots,[G:G_k]$
 are smaller than $e^{\gamma}M\log M+O(M\log\log M)$.

{\rm (iii)} The number of distinct prime divisors of
$[G:G_1],\ldots,[G:G_k]$ does not exceed $e^{\gamma}M+O(M/\log M)$.

{\rm (iv)} For the least index,
$\log[G:G_1]\ls\f{e^{\gamma}}{\log 2} M\log^2 M+O(M\log M\log\log M)$.
\endproclaim

\Proof.  Let $p_*=p_1<\cs<p_r=p^*$ be all the distinct prime
divisors of $N$.
By the supposition and Lemma 2.5,
either all the $G_i$ with $[G:G_i]\gs p_r$ are subnormal in $G$ and hence
conditions (a) and (b) in Theorem 4.1 are satisfied,
or we have  condition (c) in Theorem 4.1. In light of Theorem 4.1,
$$p_r<\max\Sb 1\ls j\ls k\\p_r\mid[G:G_j]\endSb
|\{1\ls i\ls k\colon[G:G_i]=[G:G_j]\}|\prod^r_{t=1}\f{p_t}{p_t-1}.$$
So, for some $j=1,\ldots,k$ with $[G:G_j]$
divisible by $p^*=p_r$, we have
$$\aligned&|\{1\ls i\ls k\colon [G:G_i]=[G:G_j]\}|
\\>& p_r\prod^r_{t=1}\f{p_t-1}{p_t}
=p^*\prod_{p\mid N}\f{p-1}p=(p_r-1)\prod^{r-1}_{t=1}\f{p_t-1}{p_t}
\\\gs& p_{r-1}\prod^{r-1}_{t=1}\f{p_t-1}{p_t}\gs\cs\gs p_1\f{p_1-1}{p_1}=p_*-1\endaligned$$
and hence $M\gs |\{1\ls i\ls k\colon [G:G_i]=[G:G_j]\}|
\gs1+\lfloor p^*\prod_{p\mid N}(p-1)/p\rfloor\gs p_*$.

 Note that $M>p^*\prod_{p\ls p^*}(p-1)/p$.
 Let $c(M)$ be as in the proof of Lemma 4.2. By Lemma 4.2 and its proof,
we have
$$p_1<\cs<p_r=p^*<c(M)=e^{\gamma}M\log M+O(M\log\log M)$$
and
$$r\ls\pi(p_r)=\pi(p^*)\ls \pi (c(M))=e^{\gamma}M+O(M/\log M).$$

It is known that $\zeta(2):=\sum^{\infty}_{n=1}n^{-2}=\pi^2/6$.
Let $\al(M)=2+\lfloor\log_2(\zeta(2)c(M))\rfloor$.
By induction, if $0\ls x_1,\ldots,x_n\ls 1$ then
$\prod_{i=1}^n(1-x_i)\gs1-\sum_{i=1}^nx_i$.
Thus
$$\align&\prod_{p\ls c(M)}\(1-\f1{p^{\al(M)+1}}\)\gs1-\sum_{p\ls c(M)}\f1{p^{\al(M)+1}}
\\ \gs&1-\sum_{p\ls c(M)}\f1{p^2}\cdot\f1{2^{\al(M)-1}}\gs1-\f{\zeta(2)}{2^{\al(M)-1}}
>1-\f1{c(M)}\endalign$$
and hence
$$1-\prod_{p\ls c(M)}\(1-\f1{p^{\al(M)+1}}\)<\f1{c(M)}\ls\f1M\prod_{p\ls c(M)}\f{p-1}p.$$
Therefore
$$\align &\prod_{p\ls c(M)}\sum_{n=0}^{\infty}\f1{p^n}-\f1M
=\prod_{p\ls c(M)}\f p{p-1}-\f1M
\\<&\prod_{p\ls c(M)}\f p{p-1}\cdot\prod_{p\ls c(M)}\(1-\f1{p^{\al(M)+1}}\)
\\=&\prod_{p\ls c(M)}\(\sum^{\infty}_{n=0}\f1{p^n}-\f1{p^{\al(M)+1}}
\sum^{\infty}_{n=0}\f1{p^n}\)
=\prod_{p\ls c(M)}\sum^{\al(M)}_{n=0}\f1{p^n}.\endalign$$

If $w_{\Cal A}(x)=m$ for all $x\in G$, then $\sum_{i=1}^k[G:G_i]^{-1}=m$
by Lemma 2.2 of [Su8].
Set $S=\{[G:G_1],\ldots,[G:G_k]\}$ and let $T(M)$ be the set of positive integers
which have no prime divisors greater than $c(M)$.
Then
$$1\ls\sum^k_{i=1}\f1{[G:G_i]}\ls M\sum_{n\in S}\f1n
<M\sum\Sb n\in T(M)\\n\gs[G:G_1]\endSb\f1n$$
and thus
$$\sum\Sb n\in T(M)\\n<[G:G_1]\endSb\f1n=\sum_{n\in T(M)}\f1n
-\sum\Sb n\in T(M)\\n\gs[G:G_1]\endSb\f1n
<\prod_{p\ls c(M)}\sum^{\infty}_{n=0}\f1{p^n}-\f1M
<\prod_{p\ls c(M)}\sum^{\al(M)}_{n=0}\f1{p^n}.$$
Now it is clear that $\prod_{p\ls c(M)}p^{\al(M)}$ cannot be less than $[G:G_1]$.
So $[G:G_1]\ls
\prod_{p\ls c(M)}p^{\al(M)}$ and hence
$$\log [G:G_1]\ls\al(M)\theta(c(M))\ls
l(M):=\l(2+\log_2(\zeta(2)c(M))\r)\pi(c(M))\log c(M)$$
where $\theta(x)=\sum_{p\ls x}\log p\ (\ls\pi(x)\log x)$ is
the Chebyshev $\theta$-function.
By the proof of Lemma 4.2,
$$\log(\zeta(2)c(M))=\log\f{\pi^2}6+\log M+O(\log\log M)=\log M+O(\log\log M)$$
and
$$\align\pi(c(M))\log c(M)=&c(M)+O(c(M)/\log c(M))
\\=&e^{\gamma}M\log M+O(M\log\log M)+O(e^{\gamma}M)
\\=&e^{\gamma}M\log M+O(M\log\log M).\endalign$$
So we finally have
$$\align l(M)&=\f1{\log 2}(\log M+O(\log\log M))(e^{\gamma}M\log M+O(M\log\log M))
\\&=\f{e^{\gamma}}{\log 2}M\log^2M+O(M\log M\log\log M).\endalign$$
This concludes our proof. \qed

\Remark\ 4.3. Obviously Theorem 4.3 provides more detailed information
than Theorem 1.1 does.
\medskip

For a nontrivial uniform cover (1.3) of $\Z$,
it is known that among the $k$ moduli
the largest $n_k$ occurs at least $p$ times
where $p$ is the smallest prime divisor of $n_k$
(cf. [Ne], [NZ], [Su3]). This, together with Theorem 4.3(i), suggests
the following conjecture.

\proclaim{Conjecture 4.1} Let $(1.1)$ be a nontrivial uniform
cover of a group $G$ by left cosets of subnormal subgroups.
Set $n=\max_{1\ls i\ls k}[G:G_i]$. Then
$|\{1\ls i\ls k\colon [G:G_i]=n\}|$
is not less than the least prime divisor of $n$.
\endproclaim

\enddocument